\documentclass[a4paper,11pt,reqno]{amsart}

\usepackage{a4wide}
\usepackage{latexsym}
\usepackage{amsopn,amscd}
\usepackage{amsfonts}
\usepackage{amssymb}
\usepackage{amsthm}
\usepackage{amsmath}
\usepackage{todonotes}

\usepackage[all]{xy}
\usepackage{scalefnt}
\usepackage{hyperref}

\def\can{\mathrm {can}}
\def\Aut{\mathrm{Aut}}
\def\Out{\mathrm{Out}}
\def\InAut{\mathrm{InAut}}

\long
\def\MSC#1\EndMSC{\def\arg{#1}\ifx\arg\empty\relax\else
      {\par\narrower\noindent
      2020 Mathematics Subject Classification. #1\par}\fi}

\long
\def\KEY#1\EndKEY{\def\arg{#1}\ifx\arg\empty\relax\else
    {\par\narrower\noindent
      Keywords and Phrases: #1\par}\fi}

\title[History of Lie brackets, crossed modules, Lie-Rinehart algebras]
{On the history of Lie brackets, crossed modules, and Lie-Rinehart algebras}

\author{Johannes Huebschmann}
\address{
\newline
Universit\'e de Lille - Sciences et Technologies 
\\
D\'epartement de Math\'ematiques\\
\newline CNRS-UMR 8524,
Labex CEMPI (ANR-11-LABX-0007-01)
\\
\newline
59655 Villeneuve d'Ascq Cedex, France\\
\newline
Johannes.Huebschmann@univ-lille.fr
 }

\date{\today}

\numberwithin{equation}{section}
\usepackage[capitalize]{cleveref}
\begin{document}
\setcounter{page}{1}
 
\maketitle
\centerline{Dedicated to the memory of Kirill Mackenzie}
\begin{abstract}
\noindent
The aim here is to sketch  the development of ideas
related to brackets and similar concepts:
Some purely group theoretical combinatorics due to Ph. Hall led to a proof
of the Jacobi identity for the Whitehead product in homotopy theory.
Whitehead
introduced crossed modules 
to characterize a second relative homotopy group;
guided by combinatorial group theory considerations,
Reidemeister and Peiffer explored this kind of structure
to develop normal forms for the decomposition of a 3-manifold;
but crossed modules are also 
lurking behind a forgotten approach of Turing
to the extension problem for groups: Turing concocted
the obstruction 3-cocycle
isolated later by Eilenberg-Mac Lane
and already proved 
the Eilenberg-Mac Lane theorem to the effect
that the vanishing of the class
of that cocycle
is equivalent to the existence of a solution for the
corresponding extension problem.
This Turing cocycle is  related to what has come to be
known as Teichm\"uller cocycle.
There was a parallel development for Lie algebras
including a forgotten paper by Goldberg
and, likewise, for Lie-Rinehart algebras
and Lie algebroids.
Versions of Turing's theorem were discovered 
several times
under such circumstances, 
and there is rarely a hint at 
the mutual relationship.
Also, Lie-Rinehart algebras have for long occurred in the literature
on differential algebra,
at least implicitly.

\end{abstract}

\MSC 

\noindent
Primary: 
01A60, 17B66, 18G45 

\noindent
Secondary: 
12G05, 
12H05, 
17-03, 
17B55, 
20J05, 
53-03, 
53D17, 
55Q15, 
58H05. 

\EndMSC

\KEY
Graded Lie algebra, crossed module, 
obstruction to nonabelian extension, nonabelian extension of groups,
nonabelian extension of Lie algebras, 
Lie-Rinehart algebra, Lie algebroid, differential algebra, Whitehead product,
identities among relations, Teichmueller cocycle, Galois cohomology
\EndKEY

\section{Introduction}
Lie brackets
and variants thereof have a history.
Here I attempt to give an account  
of the evolution of some ideas
intricately relating
groups and Lie algebras.
The story, not complete,
 somewhat reflects the accursed first half of the 20th century, however.
The reader can
imagine what the outcome of interaction or
collaboration among, e.g., Baer, Reidemeister, Turing, Whitehead,
Teichm\"uller, Hochschild, Magnus,  Witt, etc. might have been!

Section \ref{grlie} recalls how some combinatorics
due to P. Hall \cite{zbMATH03010343}, at first purely group 
theoretical,
eventually led to a  proof of the graded Jacobi identity
of the Whitehead product.

Section \ref{crossed} discusses 
how the
crossed module concept evolved
out of considerations of Turing, Whitehead, Reidemeister and Peiffer.
In particular,
Turing's forgotten 
paper \cite{MR1557005} 
about the extension problem for groups
{\em implicitly contains 
the idea of a (free) crossed module.
Furthermore, 
it develops the 
group cohomology $3$-cocycle which encapsulates
the obstruction for an abstract kernel to 
be extendible\/} (i.e., realizable by a group extension, see Section \ref{crossed}), and   
\cite[Theorem 4 p.~365]{MR1557005}
{\em  establishes the fact that the 
vanishing of the cohomology class is 
equivalent to extendibility of the abstract 
kernel under discussion\/}.
In the rest of the paper
I refer to this theorem as \lq\lq Turing's theorem\rq\rq.
That obstruction 
3-cocycle was later identified in \cite {MR0020996}
together with a proof of Turing's theorem.
From Section \ref{crossed} on,
Turing's theorem is a kind of Ariane's thread:
I illustrate how variants of this theorem were discovered several times
in the literature, mostly without notice of the mutual relationship.
To this end, I  spell out
various situations 
where crossed modules
provide conceptual 
means to simplify otherwise
complicated calculations.
Also, these conceptual 
means provide new insight;
for example, from crossed modules one can built a cohomology theory
for topological groups not covereed by (global) cocycle
descriptions.

Section \ref{teich} discusses the relationship
between the Teichm\"uller cocycle and crossed modules
and the eventual outcome thereof.

Section \ref{lie} recalls the corresponding development
in the Lie algebra case.
I show in particular that the forgotten paper
\cite{MR57853}
implicitly develops a Lie algebra variant of Turing's theorem
in terms of crossed 2-fold extensions of Lie algebras.

Section \ref{lr} illustrates how Lie-Rinehart algebras,
a far reaching  algebraic generalization of Lie algebroids,
have for long occurred in the literature 
on differential algebra, at least implicitly.

A MR or Zentralblatt search does not fully reveal 
historical details
of the kind  described below
nor does Google scholar,
and it proved indispensable 
to look up the papers themselves 
and to try to assimilate the content;
indeed at times those bibliographical tools
do not even find relevant references
quoted in the papers themselves.
A bibliographic metric fails to unveil the origins of the mathematical ideas
I describe below; see the penultimate paragraph of Section 
\ref{crossed} for details.
A historical irony is the fact that
Turing and Teichm\"uller, both forerunners
regarding the third group cohomology group,
worked as WWII codebreakers on opposite sides, cf.
\cite {MR1749515}, 
\cite[p.~12]{MR1152479}.\footnote{Turing, Hilton and J.H.C. Whitehead
worked as WWII codebreakers at Bletchley Park (UK)
\cite {MR1749515}. }
How the second group cohomology group
evolved out of work of Brauer \cite{zbMATH03004460}, Noether 
\cite{zbMATH02572834, zbMATH02562955}
and Baer \cite{zbMATH03012737} 
is not an issue in this paper, and the reader should not confuse
crossed modules and crossed products.

I hope my considerations 
help some of the presently young 
mathematicians a bit  to maintain contact with the past.
Various of the items I discuss in this paper
were among K. Mackenzie's research interests.

\section{Graded Lie algebras in the early days of algebraic 
topology: From
P. Hall's \lq\lq collecting process\rq\rq\ to
the Jacobi identity for the Whitehead product}
\label{grlie}
The paper \cite{MR1581549}
starts with the sentence:
\lq\lq P. Hall 
\cite{zbMATH03010343}
established a wide system of relations among higher commutators.\rq\rq\
\cite[Theorem 4.1]{MR0038336} says: {\em The standard 
commutators of weight $n$ are a basis of $F_n/F_{n+1}$.}
(Here $F_n$ is the group generated by commutators of weight $n$
and higher in a free grouyp $F$.) Thereafter:
\lq\lq ... the standard commutators are precisely those which arise
in Philip Hall's collecting
process given in \cite{zbMATH03010343}.\rq\rq

In \cite{whitethr},
J.H.C. Whitehead introduced 
the operation of
what later has come to be known as {\em Whitehead product\/},
an operation of the kind
$\pi_q \times \pi_r \to \pi_{q+r-1}$
involving the homotopy groups $\pi_*$ of a space.
For $* \geq 2$,
this operation yields 
a graded Lie algebra (the homotopy groups being regraded down by $1$) and,
for $r \geq 0$, the operation $\pi_1 \times \pi_r \to \pi_r$
gives the classical action of $\pi_1$ on $\pi_*$ ($\pi_0$ being only a set).
In his MR review of Samelson's paper \cite{MR60819},
P. Hilton wrote: \lq\lq The object of this note is to establish a formula
which is an important step towards proving the conjectured
Jacobi identity for Whitehead products.\rq\rq\ 
Thus in the mid 1950s  the idea of a graded Lie algebra was well 
understood among topologists.
The papers  \cite{MR79764}, \cite{MR68218}, \cite{MR65926}, \cite{MR64396}, 
\cite{MR0091473}, \cite{MR65927}
offer proofs of the Jacobi identity for the Whitehead product.
In his MR review on \cite{MR64396},
J. Moore wrote
\lq\lq This identity has been proved recently by others including Hilton and Serre, Massey and Uehara,
G.W. Whitehead, Nakaoka, and Toda.\rq\rq\ 
A footnote to \cite{MR0091473} says:
\lq\lq The authors have been informed that independent proofs of 
this Jacobi identity
have recently been found by H. Toda and M. Nakaoka in Japan,
by Serre, Hilton, and Green in Europe, and by G. Whitehead in this country.\rq\rq\ 
I am indebted to Serre for having 
helped me to clarify his unpublished contribution
and that of Green (\lq\lq Sandy Green\rq\rq):
At the end of the introduction of
\cite{MR68218}, Hilton wrote:
\lq\lq The author wishes to acknowledge the decisive contributions made by
J.-P. Serre and J. A. Green; the fundamental idea in the proof of Theorem A
(identifying the homotopy group $\pi_n$ ($n \geq2$) of a one point union of
finitely many spheres with a  direct sum of $n$th
homotopy groups of spheres arising from systematically exploiting
a family of Whitehead products) 
is due to Serre, and the necessary
extension of the algebraic methods is due 
essentially to Green.\rq\rq\ 
In a letter to Hilton, 
Serre had 
explained the idea that the homotopy groups of
a one point union of
finitely many spheres 
recover all multivariable operations 
among homotopy groups of spheres.
Using the graded analogue of Witt's result
\cite{MR1581553}
establishing the link between 
free Lie algebras and free associative algebras
(in modern terminology:
the universal algebra associated with a free Lie algebra is
the free associative algebra
on the very same generators),
Hilton  managed to render Serre's idea rigorous.
To this end, he introduced 
\lq\lq basic commutators\rq\rq\ 
imitating 
 Philipp Hall's basic commutators in group theory; 
these are the standard commutators in \cite{MR0038336}),
see also \cite[Theorem 5.13. A p. 343]{magkarso} and the literature there.
The corresponding somewhat more ring theoretic analogue
is \cite[Theorem 5.8 p. 323]{magkarso}, first proved in 
\cite{MR1581549}, a paper in turn quoted by Hilton.
It is presumably at this point where
Green worked out \lq\lq the necessary
extension of the algebraic method\rq\rq:
In a footnote on p.~164
of \cite{MR68218}, Hilton wrote: 
\lq\lq The ensuing calculation is implicit in the last remark of
\cite{MR60819}.\rq\rq\ 
Thus, commutator calculus known from group theory,
in particular from \cite{zbMATH03010343},
applied to the loop space, close to a group,
indeed, an H-space, properly interpreted,
apparently enabled Hilton
to prove his Theorems A and B (Jacobi identity
for Whitehead products).
In a footnote to \cite{MR68218}, 
Hilton noted that Hurewicz established the Jacobi 
identity under discussion as well, but I could not find any trace in the 
literature to this effect.
Also, in \cite{MR1581553}, Witt acknowledged that a talk of Magnus about the
results in
\cite{MR1581549} prompted Witt to write his paper.

As for the terminology {\em graded Lie algebra\/}:
In \cite{zbMATH03188566} (seminar delivered in May 1955,
typed manuscript dated \lq\lq Juillet 1957\rq\rq),
Cartier considered the notion of a graded Lie algebra
({\em alg\`ebre de Lie gradu\'ee\/}).
Actually, Cartier explored  a differential graded Lie algebra
without saying so, 
the universal algebra associated with
 that differential graded Lie algebra
yields a variant of the description 
in
\cite[Ch. XIII, Ex. 14, p.~287]{MR1731415}
of the standard resolution
of the ground field nowadays
familiar in Lie algebra theory,
 and the map $\varphi$ on the bottom of
p. 287 of \cite{MR1731415}
renders the two descriptions under discussion explicit.
For a given Lie algebra $\mathfrak g$,
that variant of the standard resolution
of the ground field 
is written in \cite[Ch. XIII, Ex. 14, p.~287]{MR1731415}
as $W(\mathfrak g)$---coincidentally, in \cite{MR0042427, MR0042426},
H. Cartan used the letter $W$ 
for the Weil algebra.
 Among other things, \cite[Ch. XIII, Ex. 15, p.~287]{MR1731415}
describes the structure of 
$W(\mathfrak g)$ as a \lq\lq graded differential algebra\rq\rq.
Thus Cartan and Eilenberg had understood, at least implicitly,
the idea of a differential graded Lie algebra in September 1953,
date of the preface of their book \cite{MR1731415}.
This book is quoted in  \cite{zbMATH03188566}.
In his MR review on \cite{MR85241},
J. Moore wrote
\lq\lq If one defines ..., one obtains a ring satisfying the usual 
identities
for a graded Lie ring.\rq\rq\ 
In a footnote to \cite{MR85241}, Hilton wrote:
\lq\lq Cartier has considered the notion of a graded Lie algebra in
\cite{zbMATH03188566}.\rq\rq\ 
Thus at the time,
the terminology  {\em graded Lie ring\/}
and
{\em graded Lie algebra\/}
was lingua franca in algebraic topology.

\section{From Schreier's extension theory to crossed modules}
\label{crossed}

To explore group extensions, Schreier
 introduced factor sets and studied their properties \cite{zbMATH02585724};
see
\cite[\S 48]{MR0266978} for an account of Schreier's theory.
Baer \cite{zbMATH03012737} elaborated on 
\cite{zbMATH02585724}
in terms of an additional ingredient,
the idea of an abstract kernel:  
For a group $Q$, an {\em abstract $Q$-kernel\/} 
or simply {\em $Q$-kernel\/} 
(terminology due to \cite {MR0020996})
is a group $N$ together with a homorphism
$Q \to \Out(N) = \Aut(N)/ \InAut(N)$
(\lq\lq Kollektivcharacter\rq\rq\ in the terminology of
\cite{zbMATH03012737}).
An extension $N \rightarrowtail E \twoheadrightarrow Q$
determines, via conjugation in $E$, an abstract $Q$-kernel
structure $Q \to \Out(N)$. However,
not every abstract kernel is {\em realizable\/}
in this manner;
already \cite{zbMATH03012737} exhibits a counterexample.

To understand a 
relative version of the product operation discussed in the previous section,
 Whitehead isolated the concept
of a crossed module;
the idea occurs in
\cite{whitethr}, is made precise in 
\cite{whitefou}, and the crossed module terminology
occurs in \cite[\S 2 p.~453]{whitefiv};
the paper \cite{whitethr}
quotes 
\cite{zbMATH02553653, zbMATH03013614}
and
\cite{whitefiv}
quotes 
\cite{zbMATH03030739}.
Here
is the definition:

A {\em crossed module\/}
consists of groups $G$ and $\Gamma$,
an action of $G$ on $\Gamma$ from the left which we here write
as 
\begin{equation*}
G \times \Gamma \to \Gamma,
\ (x,y) \mapsto {}^x y,
\end{equation*}
and a $G$-homomorphism
$\partial \colon \Gamma \to G$, 
the action of $G$ on itself being by conjugation,
subject to
\begin{equation}
 ( {}^{\partial(x)}y)xy^{-1}x^{-1} =1,\ x, y \in \Gamma.
\label{peiff}
\end{equation}

The paper \cite{MR0033287} also introduces
crossed modules (without a name),
 hypothesis (i) in I.2 (p.~738) of this paper 
encapsulates the crossed module axiom \eqref{peiff},
and Theorem 3.1 (p.~742)
establishes the equivalence of
abstract kernels and crossed modules.
Likewise,
Footnote 2 of \cite{whitefiv} says:
\lq\lq Anne Cobbe has pointed out to me that a crossed 
$(\gamma,d)$-module determines a $Q$-kernel ... and that any 
$Q$-kernel has a representation
as a crossed module.\rq\rq\

The paper \cite{peiffone} 
(Ph.D. thesis supervised by Reidemeister,
submitted for publication in 1944)
arose out of a combinatorial study of 3-manifolds and
also isolated
the identities \eqref{peiff}
defining a crossed module;
these identities have since \cite{MR47046}  come to be known
as {\em Peiffer identities\/}.
See \cite{reideone} for a survey.
But there is no hint of Whitehead's crossed modules,
nor had Whitehead recognized the Peiffer-Reidemeister theory.
In the notation of \cite[\S 1]{reideone},
$\mathfrak R/\mathfrak P \to \mathfrak S$
is a crossed module
(here $\mathfrak P$ 
refers to the corresponding Peiffer elements) and, furthermore,
in \cite[\S 2]{reideone},
$\mathfrak R/\mathfrak P \to \mathfrak S$ is the free crossed module
associated with the presentation of the group $\mathfrak G$ under discussion.
A footnote on p.~ 379 of
\cite{MR48461}
hints at the relationship and 
\cite[Section 2,  penultimate paragraph on p.~207] {MR154283}
explicitly recognizes 
the equivalence of the two approaches. 
In response to a letter,
 R. Peiffer  revealed
she was not familiar with the notion of crossed module, however.
There is no mention of Reidemeister or Peiffer
in \cite{MR405434}, no mention of Whitehead in
\cite{MR564434}. 
At the end of the 1970s, the relationship became more widely known
 \cite{MR564432, MR524183}.

In \cite[p.~423]{whitethr}, in the course of the proof that
the second relative homotopy group
arising from attaching 2-cells is a free crossed module,
i.e., characterized merely by the identities \eqref{peiff},
Whitehead gives credit to  \cite{zbMATH02587522}: A footnote
says: \lq\lq This fact is, so to speak, half the content of the proof.\rq\rq\ 
Here Whitehead  
interprets  the identities \eqref{peiff}
in terms of the Wirtinger relations
defining the fundamental group of the exterior of a 
(tame) knot in the 3-sphere.
The paper \cite{browwhit} offers a modern account of these ideas.
An extension thereof is in \cite{braids}; 
the main result of this paper
says that, for $n \geq 3$, Artin's braid group $B_n$,
as a crossed module over itself, has a single generator, which can be 
taken to be any of the Artin generators, and that the kernel
of the surjection from the corresponding free crossed module
is the second homology group $\mathrm H_2(B_n)$, well known to be
cyclic of order 2 when $n \geq 4$ and trivial for $n=2$ and $n=3$.
This includes an interpretation of the Artin relations
in terms of Peiffer identities.

To delve briefly  
into \cite{MR1557005} and fulfill the promise  in the introduction,
sticking to  Turing's notation \cite{MR1557005}
(to some extent the same as that in  \cite{zbMATH03012737}),
let $\langle e_1,\ldots, e_n; r_1,\ldots, r_l\rangle$
be a presentation of a group 
$\mathfrak G$.\footnote{Turing's notation is $\mathfrak G'$
rather than $\mathfrak G$ but, in modern notation, this is misleading
since $\mathfrak G'$ is standard notation for $[\mathfrak G,\mathfrak G]$.}
Turing gives  credit to  \cite{zbMATH02587522}
for the following requisite combinatorial group theory:
Let $\mathfrak F$ denote the free group on $e_1,\ldots, e_n$ and
$\mathfrak R$ the normal closure of 
$r_1,\ldots, r_l$ in $\mathfrak F$.
Let $\varPhi$ be the free $\mathfrak F$-operator group
on $r_1,\ldots, r_l$ and let $\mathbf P$ denote the kernel of the obvious
epimorphism $\tau\colon \varPhi \to \mathfrak R$. Turing refers to 
the members of $\mathbf P$ as \lq\lq relations between 
relations\rq\rq\ 
(in his Zbl review of Turing's paper, 
Baer attributes the idea of exploring such
relations between relations
to \cite{zbMATH02587522}); these are the identities among relations
in \cite{peiffone} and \cite{reideone},
see  \cite{MR662431}
for a modern account (but at the time of writing we were unaware of
Turing's contributions).
Write the action as
\begin{equation}
\varPhi \times \mathfrak F \longrightarrow \varPhi,\ (y,w) \mapsto y^w, 
\ y\in \varPhi,\ w\in \mathfrak F,
\end{equation}
so that
$\tau (y^w)= w^{-1}\tau (y)w$, for $y\in \varPhi$ 
and $w \in \mathfrak F$.
Accordingly, 
write the elements \cite[(9) p.~363]{MR1557005} of $\mathbf P$
in the form
\begin{equation}
y^{\tau(x)} x^{-1} y^{-1} x,
\ \ 
y= r_i^a \in \varPhi,\ \ x=  r_j^b\in \varPhi,\ \ a, b \in \mathfrak F,
\ i,j = 1, \ldots, l.
\label{turing}
\end{equation}
With the requisite modifications,  since we are now working with
a right action,
these elements formally recover those on the left-hand side of 
\eqref{peiff}.
The notation in 
\cite{MR1557005} is $E_{i,a}$ for $r_i^a \in \varPhi$ ($a \in \mathfrak F$), 
and its is useful to
remember that, as a group, $\varPhi$ is freely generated by
the $E_{i,a}$ ($i = 1, \ldots, l$), as $a \in \mathfrak F$.

Now Turing proceeds as follows
to construct a family of members of $\mathbf P$
which, together with \eqref{turing}, generate $\mathbf P$:
Choose a \lq\lq function\rq\rq\  $v \colon \mathfrak F \to \mathfrak F$
that factors through the canonical epimorphism 
$\can\colon \mathfrak F \to \mathfrak G$
as  $\mathfrak F \to \mathfrak G \to \mathfrak F$
and whose composite with the canonical epimorphism 
coincides with that epimorphism.
Thus $v$ encapsulates a section,
not necessarily a homomorphism, for that epimorphism.
For $a \in \mathfrak F$, write $v_a = v(a)$,
so that $r_a= v_a^{-1}a \in \mathfrak R$;
then, for $b \in \mathfrak F$, necessarily
$r_{v_br_i}= v_{v_br_i}^{-1}v_br_i =r_i$ ($i = 1, \ldots, l$), since
$v_{v_br_i} =v_b$.
By means of a recursive procedure,
Turing then defines corresponding members
$R_{v_br_i}$ of $\varPhi$, and 
\cite[Theorem 3 p.~364]{MR1557005} says the following:
{\em The members {\rm \eqref{turing}}
and
\begin{equation}
(R_{v_br_i} E_{i,e}^{-1})^a,\ a,b \in \mathfrak F,\ i = 1, \ldots, l,
\label{significant}
\end{equation}
generate the kernel $\mathbf P$ of $\tau \colon \varPhi \to \mathfrak R$.}
Turing applies these constructions to an abstract kernel
$(\mathfrak N, X\colon \mathfrak G \to \Out(\mathfrak N)$)
(the notation $X$ is due to
\cite{zbMATH03012737}),
together with a lift
$\chi \colon \mathfrak F \to \Aut(\mathfrak N)$
of $X$. Theorem 4 in
\cite[ p.~365]{MR1557005} says
the following 
(in somewhat more modern terminology):
{\em There is a group extension 
$\mathfrak N \rightarrowtail \mathfrak E \twoheadrightarrow \mathfrak G$
realizing the 
abstract kernel
$(\mathfrak N, X\colon \mathfrak G \to \Out(\mathfrak N))$ if and only if 
the pair $(\chi,X)$ extends to a commutative diagram
\begin{equation}
\begin{gathered}
\xymatrix
{
\varPhi \ar[r]^\tau \ar[d]_{\vartheta}
&
\mathfrak F \ar[d]^{\chi}
\ar[r]^\can
&
\mathfrak G \ar[d]^{X} \ar[r]
&
1
\\
\mathfrak N
\ar[r]_{\can\phantom{aaa}}
&
\Aut(\mathfrak N)
\ar[r]_{\phantom{}\can}
&
\Out(\mathfrak N)
\ar[r]
&
1
}
\end{gathered}
\label{turing2}
\end{equation}
with exact rows 
such that $\vartheta$ is trivial on the elements {\rm\eqref{significant}}.
}
Here \lq\lq $\can$\rq\rq\  refers to the canonical maps, and
Turing's notation is $\mathfrak G$ for $\mathfrak E$ and 
$\mathfrak G'$ for $\mathfrak G$.
Needless to point out, neither the idea of an  exact sequence
nor that of a commutative diagram was available to Turing.

Let $T$ (here $T$ stands for \lq\lq Turing\rq\rq) 
denote the normal closure in $\varPhi$ of the elements
\eqref{turing}---the various groups $T$, $\varPhi$, $\mathfrak F$
correspond precisely to the groups,
respectively,
$\mathfrak P$, $\mathfrak R$, $\mathfrak S$,
in \cite[\S 2]{reideone}.
The images of the members \eqref{significant}
generate the kernel $\pi_2$ of the induced epimorphism
$\varPhi/T \to \mathfrak R$.
Here
the notation $\pi_2$ serves here as a mnemonic
for the fact that this kernel recovers the second homotopy group
of the geometric realization  of the presentation
 (2-dimensional CW-complex realizing the presentation) 
under discussion, see, e.g., \cite{MR662431}.
With hindsight we understand that the quotient group
$\varPhi/T$
is the {\em free crossed $\mathfrak F$-module\/}
on $r_1,\ldots,r_l$ and that,
by virtue of the 
interpretation
of the third group cohomology group
in terms of \lq\lq crossed 2-fold extensions\rq\rq (details below),
 the restriction to $\pi_2$
of the  homomorphism
$\widetilde \vartheta \colon \varPhi/T \to \mathfrak N$
which $\vartheta$ induces
recovers a version of
the Eilenberg-Mac Lane 3-cocycle
 \cite {MR0020996}
associated with the abstract kernel under discussion.
Thus we can interpret 
Turing's theorem
\cite[Theorem 4 p.~365]{MR1557005}
as saying that
 {\em an abstract kernel 
is extendible if and only if
the associated $3$-cohomology class vanishes\/},
a result established nine years later
as
\cite[Theorem 8.1 p.~33] {MR0020996}.

The paper \cite{historicalnote} discusses the history of the interpretation,
for $n\geq 1$,
of the group cohomology group
$\mathrm H^{n+1}(Q,M)$ for a group $Q$ and a $Q$-module
$M$
in terms of crossed $n$-fold extensions
\cite{MR0207793, MR549932, thesis, MR576608,  
MR511461, MR2627109}; 
see also the appendix.
Suffice it to point out the following:
A {\em crossed $2$-fold extension\/}
of a group $Q$ by a $Q$-module $\pi$ is an exact sequence
\begin{equation}
0
\longrightarrow
\pi
\longrightarrow
C
\stackrel{\partial}\longrightarrow
G
\longrightarrow
Q
\longrightarrow
1
\end{equation} 
of groups
having 
$\partial \colon C \to G$
a crossed module
and
$\pi \to C$ a morphism of $G$-groups,
the action of $G$ on $\pi$ being through 
$G \to Q$.
Congruence classes of
crossed $2$-fold extensions
of $Q$ by $\pi$ 
constitute the cohomology group $\mathrm H^3(Q, \pi)$.
To reconcile this with 
Turing's theorem,
consider the pullback group $G$
\footnote{in the terminology of
\cite[p.~378 ff]{zbMATH03012737}
\lq\lq Aufl\"osung des Kollektivcharakters\rq\rq}
characterized 
by requiring that
\begin{equation}
\begin{gathered}
\xymatrix
{
G \ar[r] \ar[d]
&
\mathfrak G \ar[d]
\\
\Aut(\mathfrak N)
\ar[r]
&
\Out(\mathfrak N)
}
\end{gathered}
\end{equation}
be a pullback diagram, and let $Z$
be the center of $\mathfrak N$.
With the obvious structure,
$\mathfrak N \to \Aut(\mathfrak N)$
is a crossed module
whence it is immediate that
the homomorphism $\vartheta$ in \eqref{turing2}
 above is trivial on $T$.
Thus
the diagram \eqref{turing2} induces a morphism
\begin{equation}
\begin{gathered}
\xymatrix
{
0
\ar[r]
&
\pi_2
\ar[r]\ar[d]
&
\varPhi/T \ar[r]^{\widetilde \tau} \ar[d]_{\widetilde \vartheta}
&
\mathfrak F \ar[d]^{\widetilde \chi}
\ar[r]^\can
&
\mathfrak G \ar@{=}[d] \ar[r]
&
1
\\
0
\ar[r]
&
Z
\ar[r]
&
\mathfrak N
\ar[r] 
&
G
\ar[r] 
&
\mathfrak G 
\ar[r]
&
1
}
\end{gathered}
\label{turing3}
\end{equation}
of crossed 2-fold extensions, and the restriction
of $\widetilde \vartheta$ to $\pi_2$
recovers a version of the Eilenberg-Mac Lane 3-cocycle
\cite {MR0020996} of $\mathfrak G$
with values in the center $Z$ of $\mathfrak N$.

The sequence  \cite[(2.1) p.~738]{MR0033287}
is a crossed 2-fold extension.
The terminology \lq\lq crossed sequence\rq\rq\ is in
\cite{maclwhit}.
The interpretation of
the third group cohomology group
in terms of crossed 2-fold extensions
goes back to 
\cite{MR0207793}, without usage of the crossed module terminology,
however:
In this paper, a {\em category of interest\/}
is defined to be a subcategory $\mathfrak C$ of 
the category of all rings (associative or not)
closed under kernels, cokernels and fibered products, and an introductory 
remark says:
\lq\lq Had we considered categories of interest inside
{\em almost abelian categories\/}
(in the sense of J. Moore), then what follows could also have included 
the cases of groups and schemes.\rq\rq\ 
Now, consider two objects $A$ and $C$ of
$\mathfrak C$. An $A$-{\em structure on\/} $C$
is a split extension $C \rightarrowtail B \twoheadrightarrow A$
together with a choice of splitting.
For an $A$-structure on $C$, 
in the terminology of \cite{MR0207793},
a morphism
$\varphi \colon C \to A$ is {\em conform\/}
if the $C$-structure on $C$ which $\varphi$ induces 
coincides with canonical $C$-structure of $C$ on itself
\cite[Section 1 p.~3]{MR0207793}.

To fulfill the need for a
{\em cohomology theory that is
 satisfactory for deformation theory\/}
(\cite[Section 1 p.~3]{MR0207793}),
for objects $A$ and $M$ of $\mathfrak C$
with an $A$-module structure on $M$,
Gerstenhaber then introduces
the cohomology group
$\mathfrak E^3_{\mathfrak C}(A,M)$
in terms of congruence classes of {\em admissible\/} exact sequences 
in $\mathfrak C$
of the kind
\begin{equation}
0
\longrightarrow
M
\longrightarrow
N
\longrightarrow
B
\longrightarrow
A
\longrightarrow
0;
\end{equation}
here \lq\lq admissible\rq\rq\ signifies that
$N$ carries a $B$-structure, that $N\to B$ conforms,
and that $M \to N$ is a morphism of $B$-structures, that on $M$ being through 
$B
\to
A$
via 
the 
$A$-structure on $M$. Theorem 6 (3.) in
\cite[p.~7]{MR0207793}
includes (among other valuable insight not relevant here except 
for Lie algebras, see Section \ref{lie} below), after
the circumspect hint
\lq\lq had we given the details, the category of groups\rq\rq\ 
(i.e., $\mathfrak C$),
a natural isomorphism
$\mathfrak E^3_{\mathfrak C}(A,M) \to \mathrm H^3_{\mathfrak C}(A,M)$.
When $\mathfrak C$ is the category of groups,
the conformity constraint comes down to
\eqref{peiff} and
$\mathfrak E^3_{\mathfrak C}(A,M)$ recovers,
for a group $A$ and an $A$-module $M$,
the interpretation
of $\mathrm H^3(A,M)$
in terms of crossed 2-fold extensions.
In this case, 
\cite[Theorem 4  p.~5]{MR0207793}
is a version of
Turing's theorem.

The description of group cohomology in terms of crossed $n$-fold extensions
($n \geq 1$)
is susceptible to generalizations where
cocycles are not necessarily available.
This observation reflects the
search in
\cite[p.~3]{MR0207793} for a
cohomology theory that is
\lq\lq satisfactory for deformation theory\rq\rq.
For example, 
for an extension of a topological group $G$ by a continuous $G$-module
whose underlying bundle is non-trivial,
(global) continuous cocycles are not available,
and there are exact sequences 
relating the cocycle part of classifying in the trivial bundle case 
and a classification of the underlying bundle,
see, e.g., \cite{MR494071} and the literature there.

For a compact semisimple
 Lie group $G$, the $4$-term sequence
\begin{equation}
1
\longrightarrow
S^1 \longrightarrow \widehat \Omega G 
\longrightarrow
 \mathrm PG \longrightarrow G
\longrightarrow
1
\label{fourterm}
\end{equation}
involving the  universal central extension $\widehat \Omega G$
 of the group $\Omega G$ of based loops in $G$ by the circle group $S^1$ 
and the  group $\mathrm PG$ of based paths in $G$
acquires the structure of a crossed two-fold extension
\cite[Prop. 3.1 p.~115] {MR2366945}
and hence represents a class in a differentiable
or continuous
third cohomology group
$\mathfrak H^3(G,S^1)$
that does not admit a description in terms of global continuous 
cocycles---the corresponding cohomology being trivial 
in higher dimensions for a compact Lie group.
Also, $S^1$ is not a $G$-representation,
and an exact sequence of the kind mentioned in the previous
paragraph
does not apply.
One can view that crossed 2-fold extension as a geometric
object realizing the first Pontryagin class of the classifying space
$\mathrm BG$ or, equivalently,
the 3-cohomology class 
associated with $G$
resulting from the 3-form
which  E. Cartan exhibits in \cite[p.~197]{zbMATH02577774} 
to prove that the third Betti number of $G$ is non-zero.
With hindsight one could say
this class is an early instance of a  
3-dimensional group cohomology class.
That crossed 2-fold extension is equivalent to but 
essentially different from the corresponding gerbe, cf.
\cite{MR1197353}.
Without reference to the \lq\lq crossed\rq\rq-terminology, a 
crossed 2-fold extension occurs in
\cite[p.~ 97]{MR809126}, 
written there as 
$A
\rightarrowtail
N
\to
E
\twoheadrightarrow
G
$
together with the claim that a 3-cocycle
associated with it leads to a representation of the {}'t Hooft
commutation relations.
It seems to me
\cite{MR1197353} renders some of the claims
in \cite{MR809126} mathematically rigorous.

The crossed module concept,
variants, and generalizations thereof are  nowadays
very lively in mathematics;
see, e.g., \cite{MR2841564, MR2369164, MR1087375}
and the references there. 
The equivalence of crossed modules and 2-groups,
observed by the Grothendieck school
in the mid 1960s (unpublished), see \cite[I.1.8 p.~29]{MR2841564},
\cite{MR440553},
\cite[Theorem 5.13 p.~153]{MR2369164},
 is relevant in string theory 
\cite{MR2366945}.

Turing's contributions in \cite{MR1557005}
have been completely ignored---perhaps 
noone  except Turing himself
ever understood them.
A present day  MR citation search 
finds a single paper,
a ZBL citation search 
finds four papers but
two of them do not quote Turing's paper under discussion here
(they deal with Turing machines and cite his corresponding articles).
Google scholar finds 34 citations,
and I know a few more references
which none of these bibliographical tools find. 
I checked the references available to me and,
as far as I can see, 
Turing's paper is only cited,
without any understanding of its mathematical content.
Thus Turing would have slipped through
an evaluation system based on
bibliographic metrics.

The paper \cite{MR2316238} is the only one
which a present day  MR citation search 
finds  for  the citations of \cite{MR1557005}.
That paper has its own version of
Turing's theorem, \cite[Theorem 3.8 p.~251] {MR2316238}.
In terms of the notation of
\cite{MR2316238},
lurking behind
 \cite[Theorem 3.8 p.~251] {MR2316238}
is the commutative diagram
\begin{equation}
\begin{gathered}
\xymatrix{
&
1
\ar[r]
&\widehat N
\ar[r]
\ar@{=}[d]
&\widehat G
\ar[r]
\ar[d]
&G/N
\ar[r]
\ar@{=}[d]
&
1
\\
1
\ar[r]
&
Z
\ar[r]
&\widehat N
\ar[r]
&G
\ar[r]
&G/N
\ar[r]
&1
}
\end{gathered}
\label{neeb}
\end{equation}
encapsulating the extendibility
of the abstract kernel 
$(\widehat N, G/N \to \Out (\widehat N))$
associated with the crossed module $\widehat N \to G$
such that the cohomology class 
in $\mathrm H^3(G/N,Z)$ (suitably interpreted)
which the bottom row of \eqref{neeb}
represents is zero.
Local cocycle calculations
establish
\cite[Theorem 3.8 p.~251] {MR2316238}.
Perhaps consideration of crossed 2-fold extensions renders such calculations
obsolete. 
Besides this paper, a ZBL search also finds
\cite{MR49907}.

\section{From the Teichm\"uller cocycle to crossed modules}
\label{teich}
The paper \cite{MR0002858}
develops a group 3-cocycle
(not spelled out in this language)
over the Galois group
(associated with the data there);
the cohomology
class thereof recovers the obstruction
for the Brauer class of a central simple algebra,
normal in the sense of Galois theory,
to be trivially normal in the sense that
the Brauer class arises by extension of scalars from the fixed field.
This is another instance of  Turing's theorem discussed in 
Section \ref{crossed} 
above.\footnote{In \cite[p.~148]{MR0002858}, the precise statement reads:
\lq\lq $\xi_{\lambda,\mu,\nu}$ zerf\"allt dann und nur dann, wenn 
die Algebrenklasse
von $A$ durch Erweiterung einer Algebrenklasse \"uber $P$ entsteht.\rq\rq}
Eilenberg-Mac Lane \cite{MR0025443} rework this approach in terms of
ordinaray group cohomology.
Theorem 10.1 (p.~13) is this paper
recovers the version of Turing's theorem under discussion,
as does the exactness 
at $\mathrm H^2(K,L^*)^G$
of the exact sequence in
\cite[Section 5 p.~130]{MR0052438}.
Section 15 p.~19 ff of
 \cite{MR0025443}
also de\-velops an interpretation in terms of abstract kernels, and
the interpretation 
of Turing's theorem
in Section \ref{crossed} above
sheds new light on
 the remark
in \cite[IV.11 p.~137]{maclaboo} saying
\lq\lq The 3-dimensional cohomology groups of a group were first considered
by Teichm\"uller \cite{MR0002858}\rq\rq.

A crossed 2-fold extension 
description yields the generalized  
\lq\lq Teichm\"uller cocycle\rq\rq\  
map \cite[(8.1) p.~58]{MR3769365};
this map includes the map $\rho$ in
\cite[Theorem 3.4 (i)]{MR1803361},
given there by a cocycle description.
The exactness 
at $\mathrm{XB}(S,Q)$
of  
\cite[(9.1) p.~64]{MR3769365}
can be seen as an instance of Turing's theorem under 
the circumstances there, as can
the exactness 
at $\mathrm{XB}(T|S;G,Q)$
of the
exact sequence  \cite[(16.1) p.~104]{ MR3769366}.
The same kind of remark applies to the exact sequences
in \cite[Theorem 4.2]{MR1803361}
involving
$QB_0(R,\Gamma)$
and $QB(R,\Gamma)$, constructed
there by  cocycle consideration.
These exact sequences come down to that in
\cite[Section 5 p.~130]{MR0052438}
under the circumstances of that paper.

\section{Lie algebra crossed modules}
\label{lie}

The axiom \eqref{peiff} 
makes perfect sense for Lie algebras, and
the interpretation of the $(n+1)$th Lie algebra cohomology
group
in terms of crossed $n$-fold extensions ($n \geq 1$) is available.
Also the abstract kernel concept extends to Lie algebras
in an obvious manner, as does the equivalence
between abstract kernels and crossed modules.
The following
quote from \cite[p.~698]{MR63362}, given there as
ordinary text prose, not as a formal statement, 
spells out the Lie algebra analogue of Turing's theorem discussed 
in Section \ref{crossed}
above:
\lq\lq With every $L$-kernel $M$,
one can associate a certain 3-dimensional cohomology class for $L$ 
in the center of $M$, and the $L$-kernel is extendible
if and only if this cohomology class is $0$. This is the precise analogue 
of a result of Eilenberg-Mac Lane 
\cite {MR0020996},
for the case of groups, and can easily be proved in the same way,
mutatis mutandis.\rq\rq\ The reader will notice that
the sentence \lq\lq the $L$-kernel is extendible
if and only if this cohomology class is $0$\rq\rq\ 
is the corresponding analogue of Turing's theorem.

The paper \cite{MR57853} 
characterizes the third Lie algebra cohomology group of a Lie algebra
in terms of its nonabelian extensions
by cocycle identites which,
according to the author's claims,
are analogous to the identities in 
Teichm\"uller's theory,
cf. Section \ref{teich} above.
In the notation of \cite{MR57853}, the 
data in that paper fit into a crossed 2-fold extension
\begin{equation}
0
\longrightarrow
W
\longrightarrow
U
\longrightarrow
L^+
\longrightarrow
L
\longrightarrow
0
\label{extl}
\end{equation}
of Lie algebras,
and
\cite[Theorem 4.9 p.~476]{MR57853} 
admits the following interpretation:
There is an extension 
$
U
\rightarrowtail
L^{++}
\twoheadrightarrow
L
$
of Lie algebras
together with an epimorphism
$
L^{++}
\to
L^+
$
of Lie algebras rendering commutative the diagram
\begin{equation}
\xymatrix{
&
0
\ar[r]
&U
\ar[r]
\ar@{=}[d]
&L^{++}
\ar[r]
\ar[d]
&L
\ar[r]
\ar@{=}[d]
&
0
\\
0
\ar[r]
&
W
\ar[r]
&U
\ar[r]
&L^+
\ar[r]
&L
\ar[r]
&0
}
\end{equation}
if and only if the class in $\mathrm{H}^3(L,W)$
which \eqref{extl}
represents is zero.
The reader will notice that this is the Lie algebra version
of Turing's theorem.
Thus, contrary to
the last claim in the MR review of
\cite{MR57853},
the cocycle identities in this paper
are analogous to those in
Teichm\"uller's theory. 
The Lie algebra version
of Turing's theorem
in terms of  Lie algebra abstract kernels is in
\cite{MR63362}, with a special emphasis
on the restricted Lie algebra case, which exhibits  subtleties, as well as in
\cite[Lemma 5 p.~177]{MR57854}.

The interpretation of
the third Lie algebra cohomology group
in terms of crossed 2-fold extensions
(without usage of the crossed module terminology)
likewise goes back to 
\cite{MR0207793}:
The category of Lie algebras
is a category of interest in the sense of
\cite{MR0207793}, 
and an admissible sequence is, then,
precisely
a crossed 2-fold extension of Lie algebras.
Theorem 6 (3.) in \cite[p.~7]{MR0207793}
includes, for a Lie algebra $\mathfrak g$ and a $\mathfrak g$-module $M$, 
a natural isomorphism
$\mathfrak E^3_{\mathfrak C}(\mathfrak g,M) \to 
\mathrm H^3_{\mathfrak C}(\mathfrak g,M)$.
Over a field or, more generally,
when $\mathfrak C$ is the category 
of Lie algebras
over a commutative ring $R$ with unit having underlying $R$-module
projective,
this comes down to an isomorphism
$\mathfrak E^3_{\mathfrak C}(\mathfrak g,M) \to \mathrm H^3(\mathfrak g,M)$
onto the ordinary 3rd Lie algebra cohomology group.
(In the general case, the isomorphism
$\mathfrak E^3_{\mathfrak C}(\mathfrak g,M) \to \mathrm H^3_{\mathfrak C}(\mathfrak g,M)$ recovers the corresponding relative
Lie algebra cohomology group.)
In the Lie algebra case, 
\cite[Theorem 4  p.~5]{MR0207793}
is a Lie algebra version of
Turing's theorem.

As far as I know, \cite{MR694130}
is the first paper to use the terminology
\lq\lq Lie algebra crossed module\rq\rq.
The results in
\cite{MR641328, MR597986} extend to the Lie algebra case
in a straightforard manner.
Hochschild's wording
at the end of the first paragraph of the present section applies here.

Crossed Lie algebras occur as infinitesimal objects associated with
2-groups. 
For a semisimple Lie algebra,
a  
purely algebraic construction
(i.e., not involving a loop Lie algebra)
yields a crossed 2-fold extension representing
the infinitesimal version of
E. Cartan's 3-cohomology class mentioned in Section \ref{crossed}.
Such an extension represents the infinitesimal object
associated with the class represented by a 
crossed 2-fold extension of the kind \eqref{fourterm}, viz.
$
S^1 \rightarrowtail \widehat \Omega G 
\to
 \mathrm PG \twoheadrightarrow G$.

\section{Lie-Rinehart algebras vs Lie algebroids}
\label{lr}
K. Mackenzie studied abstract kernels of Lie algebroids
\cite[p.~221 ff]{MR896907}---his terminology 
is {\em coupling\/}---and his
Cor. 3.22 (p.~225) is
a Lie algebroid version of Turing's theorem discussed 
in Section \ref{crossed}.
In the setting of holomorphic Lie algebroids,
\cite[Theorem 3.2 (i)]{MR3345517} is the corresponding variant
of Turing's theorem.
Crossed modules of Lie algebroids
occur in 
\cite[pp.~309, 332]{MR2157566}.

Let $R$ be a commutative ring and $A$
a commutative $R$-algebra.
An $(R,A)$-{\em Lie algebra\/} \cite{rinehart}
is an $R$-Lie algebra $(L,[\,\cdot\, ,\,\cdot \,])$ together with
an action of $L$ on $A$ by derivations
and a left $A$-module structure
on $L$ that satisfy the two axioms
\begin{align}
[\alpha, a \beta] &=\alpha(a) \beta + a [\alpha, \beta],
\ a \in A,\ \alpha, \beta \in L,
\label{first}
\\
(a\alpha)b &= a (\alpha(b)),\ 
\ a,b \in A,\ \alpha \in L.
\label{second}
\end{align}
The basic example 
is the pair $(A,\mathrm{Der}(A))$,
 with the obvious structure of mutual interaction.
The concept of an $(R,A)$-Lie algebra
is the algebraic analogue of a Lie algebroid:
Indeed, a standard example is the pair
$(C^{\infty}(M), \Gamma(\lambda))$
that consists of the smooth functions
$C^{\infty}(M)$ and smooth 
sections $\Gamma(\lambda)$
of  a Lie algebroid $\lambda\colon E \to M$ 
on a smooth manifold $M$, with the obvious structure of mutual interaction.
However, there are
$(R,A)$-Lie algebras
that do not arise from a Lie algebroid.
For example, 
the $(R,A)$-Lie algebra
structure on the $A$-module of formal differentials of $A$
associated with a Poisson structure on $A$
\cite{MR1058984} does not in general arise from a Lie algebroid.
Other examples follow below.
The $R$-algebra $A$ being fixed,
$(R,A)$-Lie algebras constitute a category.
When we let the algebra variable vary, we also obtain a category.
For a good notion of morphism in this case see \cite{MR1235995}.
To cope with that situation,
in \cite{MR1058984}, I introduced the terminology
{\em Lie-Rinehart algebra\/} for 
$(R,A)$-Lie algebras
when the variable $A$ is allowed to vary.
Unfortunately, in the subsequent literature,
the two notions
 Lie-Rinehart algebra
and
$(R,A)$-Lie algebra became confused.

Lie-Rinehart algebras 
have played a major role in differential algebra
and differential Galois theory
for long, even though (it seems to me) 
the structure was not explicitly recognized.
Indeed,
a differential field $K$ with field of constants $k$
is tantamount to the $(k,K)$-Lie algebra
$(K,\mathrm{Der}_k(K))$.
Suffice it to mention the following:
Lie-Rinehart algebras occur in \cite{MR0070961} and \cite{MR10554}
 though without the name
and yield a crucial tool in
\cite{MR675306}
(the terminology there is $K|k$-Lie algebra)
for  classifying the Lie algebras that arise as
Lie algebras of differential formal groups of Ritt:
There are two cases, finite-dimensional simple 
Lie algebras and Lie algebras of Cartan type
(associated with pseudogroups of transformations---the concept of a
pseudogroup of transformations evolved out of Lie's ideas).
In the latter case, the classification is somewhat delicate
and leads to a mathematically exceedingly interesting theory
involving linear compactness.
In his MR review on \cite{MR675306}, 
A. M. Vinogradov wrote:
\lq\lq The authors begin their paper with the sentence:
\lq This is an attempt
to understand the last four papers of J. Ritt.\rq\ 
Their work is undoubtedly
an important contribution to this area. But this paper itself requires 
similar \lq attempts to understand\rq. In the reviewer's opinion, 
the use of the language of differential calculus in commutative algebras
could be very helpful for this purpose.\rq\rq\ 
The four papers under discussion are
\cite{MR1507334, MR0034759,  MR0035763, MR0037308}.
The Lie-Rinehart technology provides
some such differential calculus.

A very large number of authors
developed the idea of a Lie-Rinehart algebra,
before it was so named,
 most of whom
independently proposed their own terminology; see
\cite[p.~305]{MR2075590} (incomplete list,
compiled with the help of K. Mackenzie)
as well as \cite[p.~100] {MR1325261},
the terminology there being 
{\em Lie pseudoalgebra\/}.
This is also the terminology in
\cite {MR55323, MR55324};
these  describe the structure in a form which renders its
generality clear. 
The Lie pseudoalgebra concept in 
these two papers is not equivalent to that of a Lie-Rinehart algebra,
however, and there is no mention of the difference
 in the subsequent literature:
For a field $K$, not necessarily commutative,
\cite {MR55323} defines a $K$-Lie pseudoalgebra
(pseudo-alg\`ebre de Lie)
to be a $K$-vector space $E$
together with a Lie bracket
over the integers, 
subject to a variant 
\cite[III]{MR55323}
of
\eqref{first} with $K$ substituted for $A$,
spelled out in terms of
the operation $K \times E \to E$
of scalar multiplication.
There then results
 a $\mathbb Z$-linear map $E \to \mathrm{Der}(K)$.
For a proper skew field $K$, 
the commutator operation  turns $K$ into a 
$K$-Lie pseudoalgebra, and
an observation in \cite[Section 4] {MR55324}
says that, for a 
$K$-vector space $E$ and a
$K$-linear map
$\varphi \colon E \to K$, setting
$u(\lambda)= \varphi(u) \lambda - \lambda \varphi(u)$
($u \in E$, $\lambda \in K$) yields a $K$-Lie pseudoalgebra
in such a way that $\varphi$ is a morphism of
$K$-Lie pseudoalgebras
and that every $K$-Lie pseudoalgebra arises in this manner.
In particular, a
$K$-Lie pseudoalgebra is then an ordinary Lie algebra
over the center of $K$.
For a commutative field $K$, \cite[Section 5] {MR55324}
contains the observation that a 
 $K$-Lie pseudoalgebra $E$ of $K$-dimension at least equal to $2$
necessarily satisfies the Lie-Rinehart
axiom \eqref{second} and
is hence an ordinary $(K,k)$-Lie algebra where $k\subseteq K$
refers to the fixed field.
Nowadays there are researchers trying to develop the
Lie-Rinehart concept 
relative to  a not necessarily commutative algebra.
Also,
\cite[(1.1) p.~35]{MR0059050} 
isolates  axiom \eqref{first}
to explore connections on a Lie group and on a  homogeneous space.
I chose the terminology  {\em Lie-Rinehart algebra\/} 
in \cite{MR1058984} 
for the following reason:
The paper \cite{rinehart} goes beyond the evident formal 
similarity of the Cartan-Chevalley-Eilenberg and de Rham complexes 
noted, e.g., in \cite{MR0125867}
but presumably folk-lore at the time,
establishes the Poincar\'e-Birkhoff-Witt theorem
for these objects and thereby subsumes
ordinary Lie algebra and de Rham cohomology
under a single theory, that of derived functors.
The paper \cite{rinehart} paves, furthermore, 
the way for the subsequent interplay between
Lie-Rinehart algebras and Lie algebroids in the literature.

The Lie algebroid terminology goes back to
\cite{MR216409} and 
Lie-Rinehart algebras occur there as
Lie pseudoalgebras,
with reference to
earlier consideration of the structure.
On p.~246 the reader finds a hint that
the idea of a Lie algebroid 
underlies earlier work of
Ehresmann, Libermann, and Rodrigues,
and a remark on p.~101 of
\cite{MR1325261}
says that the Lie algebroid concept is
implicit in Ehresmann's work on higher order connections and prolongations 
of structures on manifolds at the early 1950s.

Lie-Rinehart algebras play a major role in a number of areas in 
mathematics and physics.
See \cite{MR2075590} for an overview.
Here are a few remarks:
Lie-Rinehart algebras underly a program
aimed at developing quantization
in the presence of singularities, see
\cite{MR2883413}
and the literature there.
The Schouten-Nijenhuis bracket makes sense for Lie-Rinehart algebras
\cite{MR1625610} and thereby
explains various structures that have arisen in the literature,
see also \cite{MR837203},
and so does the Fr\"olicher-Nijenhuis bracket.
The formalism of connection and curvature extends to Lie-Rinehart algebras
\cite{MR1058984}, see also
\cite{koszultw} (first published in 1965).
Differential graded Lie-Rinehart algebras
\cite{twilled, MR1764437}
underly Hodge theory,
the Dolbeault complex,
the BRST-complex \cite{MR1156541},
explain Kodaira-Spencer deformation theory, etc.

\section{My scientific relationship with K. Mackenzie}

In 1989 or so,
I established contact with K. Mackenzie
and invited him for a stay in Heidelberg.
Also we met in England at various occasions
and we intensely discussed mathematics by email.
He had a great interest in  work related to
Lie-Rinehart algebras, see, e.g.,
\cite{MR1235995}, \cite{MR1764437}, \cite{MR1325261}.
We refereed each other's work.
Once I received a request to  referee
a book project. I wrote my review, sent it in with 
a negative recommendation,
and also, at the editor's request,
suggested K. Mackenzie as a possible referee.
Later I inquired and learnt from Kirill he had refereed that project
in the first place, sent in a negative review, and suggested me
as an alternate reviewer.
The book project was realized, however.
Once I talked to Anne Kostant (then Editorial Director, Springer Mathematics) 
about this story, and here laconic reaction was:
\lq\lq The book project had been refereed.\rq\rq 

I was saddened to hear we had lost Kirill.

\section*{Acknowledgements}

I am indebted to J. Stasheff for a number of most helpful comments
on a draft of the paper.
I gratefully acknowledge support by the CNRS and by the
Labex CEMPI (ANR-11-LABX-0007-01).

\section*{Appendix}

Headings (1)-(6) below
are a few observations related to \cite{historicalnote},
(2)-(6)
complementing 
the historical remarks in \cite{MR2841564};
in none of the papers mentioned in 
(2)-(5) does the crossed module terminology occur.

\begin{enumerate}

\item
The paper \cite{MR160807} establishes a description
of the cohomology group
$\mathrm H^{n+1}(Q,M)$ ($n \geq 2$) for a group $Q$ and a $Q$-module
$M$
in terms of $n$-fold extensions of the kind
\begin{equation}
0
\longrightarrow
M
\longrightarrow
M_n
\longrightarrow
\ldots
\longrightarrow
M_2
\longrightarrow
G
\longrightarrow
Q
\longrightarrow
1
\end{equation}
having $M_2$, ..., $M_n$
ordinary $Q$-modules, that is, do not properly involve crossed modules.
Such an  $n$-fold extension 
is a crossed  $n$-fold extension
but not the most general one (and there is no mention of crossed modules 
in \cite{MR160807});
cf. the \lq\lq Note added in proof\rq\rq\ in \cite{historicalnote}.

\item
No structure equivalent to a crossed module does occur in 
\cite{MR476835}, at least not explicitly.

\item

The first satellite in
 \cite[Prop. 4.1 p.~ 313]{MR245647}
involves,  in the notation of this paper,
the crossed 2-fold extension
\begin{equation}
0
\longrightarrow
Z
\longrightarrow
E
\longrightarrow
B_0
\longrightarrow
B_1
\longrightarrow
1 .
\end{equation}
The construction of this satellite 
relies on the group case in \cite{MR0207793}.

\item
In
\cite[II.2 p.~131 ff]{MR435250}, 
\begin{equation}
0
\longrightarrow
B
\longrightarrow
K
\longrightarrow
D
\longrightarrow
G
\longrightarrow
1
\end{equation}
is a crossed 2-fold extension.
Here the authors elaborate on 
(the group case in) \cite{MR0207793, MR245647}.
The authors' aim is to study group cohomology
in a variety of groups.
Such a cohomology theory is not necessarily accessible
in terms of cocycles.
This is related to 
Gerstenhaber's
search in
\cite[p.~3]{MR0207793} for a
cohomology theory that is
\lq\lq satisfactory for deformation theory\rq\rq\ already hinted at 
in Section \ref{crossed}.

\item Section 4 p. 60 ff of
\cite{MR0393196}, 
discusses objects equivalent to
crossed $n$-fold extension;
see
in particular
p. 71 l. -8$/$-7: \lq\lq ... classifies certain kinds of $n$-fold extensions
of $X$ by the $X$-module $\Pi$.\rq\rq\ 
Such an $n$-fold extension arises
as the Moore complex of a $K(\pi,n)$-torsor.

\item 
In \cite{MR548118}
the interpretation of the third group cohomology group
in terms of crossed 2-fold extensions 
leads to a structural result in the theory of von Neumann algebras.

\item
References \cite{MR1501972, MR7108, MR0019092, MR0020996, 
 MR0266978, MR0349821, MR0033287} 
 quote
\cite{MR1557005} as does 
\cite[p.~437] {MR2841564} with the remark
\lq\lq 
the use of identities among relations for discussing nonabelian extensions
was given in \cite{MR1557005}\rq\rq.
Neither MR nor Zentralblatt find those quotes.
Google scholar finds them except \cite{MR0266978, MR0349821} 
and finds some others as well.

\item Lie-Rinehart algebras are lurking behind $D$-modules.

\item There is also a literature on crossed modules of Lie-Rinehart algebras.

\end{enumerate}

\newcommand{\etalchar}[1]{$^{#1}$}
\def\cprime{$'$} \def\cprime{$'$} \def\cprime{$'$} \def\cprime{$'$}
  \def\cprime{$'$} \def\cprime{$'$} \def\cprime{$'$} \def\cprime{$'$}
  \def\dbar{\leavevmode\hbox to 0pt{\hskip.2ex \accent"16\hss}d}
  \def\cprime{$'$} \def\cprime{$'$} \def\cprime{$'$} \def\cprime{$'$}
  \def\cprime{$'$} \def\Dbar{\leavevmode\lower.6ex\hbox to 0pt{\hskip-.23ex
  \accent"16\hss}D} \def\cftil#1{\ifmmode\setbox7\hbox{$\accent"5E#1$}\else
  \setbox7\hbox{\accent"5E#1}\penalty 10000\relax\fi\raise 1\ht7
  \hbox{\lower1.15ex\hbox to 1\wd7{\hss\accent"7E\hss}}\penalty 10000
  \hskip-1\wd7\penalty 10000\box7}
  \def\cfudot#1{\ifmmode\setbox7\hbox{$\accent"5E#1$}\else
  \setbox7\hbox{\accent"5E#1}\penalty 10000\relax\fi\raise 1\ht7
  \hbox{\raise.1ex\hbox to 1\wd7{\hss.\hss}}\penalty 10000 \hskip-1\wd7\penalty
  10000\box7} \def\polhk#1{\setbox0=\hbox{#1}{\ooalign{\hidewidth
  \lower1.5ex\hbox{`}\hidewidth\crcr\unhbox0}}}
  \def\polhk#1{\setbox0=\hbox{#1}{\ooalign{\hidewidth
  \lower1.5ex\hbox{`}\hidewidth\crcr\unhbox0}}}
  \def\polhk#1{\setbox0=\hbox{#1}{\ooalign{\hidewidth
  \lower1.5ex\hbox{`}\hidewidth\crcr\unhbox0}}}

\end{document}